\documentclass[
	5p,
]{elsarticle}

\usepackage{amsmath}

\usepackage{amssymb}
\usepackage{mathtools}
\usepackage{accents}	
\DeclarePairedDelimiter{\diagfences}{(}{)}
\newcommand{\diag}{\operatorname{diag}\diagfences}

\input{preamble}

\begin{document}

\begin{frontmatter}

\title{Towards model based control of the Vertical Gradient Freeze crystal growth process
}

\begin{abstract}
	In this contribution tracking control designs using output feedback are presented for a
	two-phase \acrlong{sp} arising in the modeling of the \acrlong{vgf} process.
	The two-phase \acrlong{sp}, consisting of two coupled \acrlongpl{fbp},
	is a vital part of many crystal growth processes due to the temporally
	varying extent of the solid and liquid 	domains during growth.
	After discussing the special needs of the process, collocated as well as
	flatness-based state feedback designs are carried out.
	To render the setup complete, an observer design is performed, using a
	flatness-based approximation of the original \gls{dps}.
	The quality of the provided approximations as well as the performance of the
	open and closed loop control setups is analysed in several simulations.
\end{abstract}


\end{frontmatter}


\newcommand{\diffinline}[2] {
        \ensuremath{\text{d}#1/\text{d}#2}
}

\newcommand{\gen}{%
	\ensuremath{\circ}
}

\newcommand{\zb}{%
	\ensuremath{\Gamma_{\gen}}
}
\newcommand{\zbs}{%
	\ensuremath{\Gamma_{\mathrm{s}}}
}
\newcommand{\zbl}{%
	\ensuremath{\Gamma_{\mathrm{l}}}
}
\newcommand{\zbt}{%
	\ensuremath{\tilde\Gamma_{\gen}(t)}
}
\newcommand{\zbts}{%
	\ensuremath{\tilde\Gamma_{\mathrm{s}}(t)}
}
\newcommand{\zbtl}{%
	\ensuremath{\tilde\Gamma_{\mathrm{l}}(t)}
}
\newcommand{\zbb}{%
	\ensuremath{\beta_{\gen}(t)}
}
\newcommand{\zbbsq}{%
	\ensuremath{\beta_{\gen}^2(t)}
}
\newcommand{\zbbs}{%
	\ensuremath{\beta_{\mathrm{s}}(t)}
}
\newcommand{\zbbl}{%
	\ensuremath{\beta_{\mathrm{l}}(t)}
}

\newcommand{\dzb}{%
	\ensuremath{\delta_{\bar z}}
}

\newcommand{\errcomplete}[3]{%
	\ensuremath{\epsilon_{\mathrm{#1}}^{#2}%
		\ifthenelse{\equal{#3}{}}{}{(#3)}%
	}
}
\newcommand{\err}[1]{
	\ensuremath{\errcomplete{#1}{}{t}}%
}

\newcommand{\zicomplete}[4][]{%
	\ensuremath{%
			\gamma%
			_{\mathrm{#2}}^{#3}%
		\ifthenelse{\equal{#4}{}}{}{(#4)}%
	}
}
\newcommand{\vicomplete}[3]{%
	\ensuremath{\dot\gamma_{\mathrm{#1}}^{#2}%
		\ifthenelse{\equal{#3}{}}{}{(#3)}%
	}
}

\newcommand{\zi}[1][t]{%
	\ensuremath{\zicomplete{}{}{#1}}
}
\newcommand{\zio}{%
	\ensuremath{\zicomplete{}{0}{}}
}
\newcommand{\vi}[1][t]{%
	\ensuremath{\vicomplete{}{}{#1}}
}
\newcommand{\vio}{%
	\ensuremath{\vicomplete{}{0}{}}
}

\newcommand{\zir}{%
	\ensuremath{\zicomplete{r}{}{t}}
}
\newcommand{\vir}{%
	\ensuremath{\vicomplete{r}{}{t}}
}
\newcommand{\vimax}{%
	\ensuremath{\vicomplete{max}{}{}}
}

\newcommand{\ziest}{%
	\ensuremath{\hat{\gamma}(t)}
}
\newcommand{\viest}{%
	\ensuremath{\dot\hat{\gamma}(t)}
}

\newcommand{\zierr}{%
	\ensuremath{\gamma_{\mathrm{e}}(t)}
}
\newcommand{\vierr}{%
	\ensuremath{\dot\gamma_{\mathrm{e}}(t)}
}
\newcommand{\zierro}{%
	\ensuremath{\gamma_{\mathrm{e}}(0)}
}
\newcommand{\vierro}{%
	\ensuremath{\dot\gamma_{\mathrm{e}}(0)}
}

\newcommand{\tm}{%
	\ensuremath{T_{\mathrm{m}}}%
}
\newcommand{\tnz}{%
	\ensuremath{T(z, t)}%
}
\newcommand{\torig}{%
	\ensuremath{T}%
}
\newcommand{\torigg}{%
	\ensuremath{T_{\gen}}%
}
\newcommand{\torigz}{%
	\ensuremath{\torig(z, t)}%
}
\newcommand{\toriggz}{%
	\ensuremath{\torigg(z, t)}%
}
\newcommand{\ts}{%
	\ensuremath{T_{\mathrm{s}}}%
}
\newcommand{\tl}{%
	\ensuremath{T_{\mathrm{l}}}%
}

\newcommand{\tref}{%
	\ensuremath{T_{\mathrm{r}}}%
}
\newcommand{\trefz}{%
	\ensuremath{\tref(z, t)}%
}
\newcommand{\trefg}{%
	\ensuremath{T_{\gen, \mathrm{r}}}%
}
\newcommand{\trefgz}{%
	\ensuremath{\trefg(z, t)}%
}
\newcommand{\trefs}{%
	\ensuremath{T_{\mathrm{s,r}}}%
}
\newcommand{\trefl}{%
	\ensuremath{T_{\mathrm{l,r}}}%
}

\newcommand{\test}{%
	\ensuremath{\hat\torig}%
}
\newcommand{\testz}{%
	\ensuremath{\test(z, t)}%
}

\newcommand{\tapp}{%
	\ensuremath{\bar{\torig}_{\gen}^N}%
}
\newcommand{\tappz}{%
	\ensuremath{\tapp(\bar z, t)}%
}
\newcommand{\tapps}{%
	\ensuremath{\bar{\torig}_{\mathrm{s}}^N}%
}
\newcommand{\tappl}{%
	\ensuremath{\bar{\torig}_{\mathrm{l}}^N}%
}

\newcommand{\zt}{%
	\ensuremath{\tilde z}
}
\newcommand{\ttrans}{%
	\ensuremath{\tilde T}%
}
\newcommand{\ttransg}{%
	\ensuremath{\tilde T_{\gen}}%
}
\newcommand{\ttranss}{%
	\ensuremath{\ttrans_\text{s}}%
}
\newcommand{\ttransl}{%
	\ensuremath{\ttrans_\text{l}}%
}
\newcommand{\ttransz}{%
	\ensuremath{\ttrans(\tilde z, t)}%
}
\newcommand{\ttransgz}{%
	\ensuremath{\ttransg(\tilde z, t)}%
}
\newcommand{\tfix}{%
	\ensuremath{\bar{\torig}_{\gen}}%
}
\newcommand{\tfixz}{%
	\ensuremath{\tfix(\bar z, t)}%
}
\newcommand{\tfixs}{%
	\ensuremath{\bar{\torig}_{\mathrm{s}}}%
}
\newcommand{\tfixl}{%
	\ensuremath{\bar{\torig}_{\mathrm{l}}}%
}

\newcommand{\dzi}{%
	\ensuremath{\Delta\zi}%
}
\newcommand{\Vdzi}{%
	\ensuremath{\Psi_0^t(\Delta\gamma)}%
}
\newcommand{\dvi}{%
	\ensuremath{\Delta\vi}%
}
\newcommand{\dai}{%
	\ensuremath{\Delta\ddot\gamma(t)}%
}
\newcommand{\treft}{%
	\ensuremath{\tref(z -\Delta\zi, t)}%
}
\newcommand{\trefgt}{%
	\ensuremath{\trefg(z -\Delta\zi, t)}%
}
\newcommand{\trefst}{%
	\ensuremath{\trefs(z -\Delta\zi, t)}%
}
\newcommand{\treflt}{%
	\ensuremath{\trefl(z -\Delta\zi, t)}%
}

\newcommand{\terrn}{%
	\ensuremath{e}%
}
\newcommand{\terr}{%
	\ensuremath{\tilde e}%
}
\newcommand{\terrz}{%
	\ensuremath{\terr(z, t)}%
}
\newcommand{\terrsq}{%
	\ensuremath{\terr^2}%
}
\newcommand{\terrzsq}{%
	\ensuremath{\terr^2(z, t)}%
}

\newcommand{\terrs}{%
	\ensuremath{\terr_{\mathrm{s}}}%
}
\newcommand{\terrssq}{%
	\ensuremath{\terrs^2}%
}
\newcommand{\terrso}{%
	\ensuremath{\terrs(0,t)}%
}
\newcommand{\terrsz}{%
	\ensuremath{\terrs(z,t)}%
}
\newcommand{\terrszsq}{%
	\ensuremath{\terrs^2(z,t)}%
}
\newcommand{\terrszbsq}{%
	\ensuremath{\terrs^2(\zbs,t)}%
}

\newcommand{\terrl}{%
	\ensuremath{\terr_{\mathrm{l}}}%
}
\newcommand{\terrlsq}{%
	\ensuremath{\terrl^2}%
}
\newcommand{\terrlo}{%
	\ensuremath{\terrl(0,t)}%
}
\newcommand{\terrlz}{%
	\ensuremath{\terrl(z,t)}%
}
\newcommand{\terrlzsq}{%
	\ensuremath{\terrl^2(z,t)}%
}
\newcommand{\terrlzbsq}{%
	\ensuremath{\terrl^2(\zbl, t)}%
}

\newcommand{\terrg}{%
	\ensuremath{\terr_{\gen}}%
}
\newcommand{\terrgz}{%
	\ensuremath{\terrg(z,t)}%
}
\newcommand{\terrgsq}{%
	\ensuremath{\terrg^2}%
}
\newcommand{\terrgzsq}{%
	\ensuremath{\terrg^2(z, t)}%
}
\newcommand{\terrgzbsq}{%
	\ensuremath{\terrg^2(\zb, t)}%
}

\newcommand{\cgain}{%
	\ensuremath{\kappa}%
}
\newcommand{\cgaing}{%
	\ensuremath{\cgain_{\gen}}%
}
\newcommand{\cgains}{%
	\ensuremath{\cgain_{\mathrm{s}}}%
}
\newcommand{\cgainl}{%
	\ensuremath{\cgain_{\mathrm{l}}}%
}

\newcommand{\hc}{%
	\ensuremath{k_{\gen}}
}
\newcommand{\hcs}{%
	\ensuremath{k_{\mathrm{s}}}
}
\newcommand{\hcl}{%
	\ensuremath{k_{\mathrm{l}}}
}

\newcommand{\hd}{%
	\ensuremath{\alpha_{\gen}}
}
\newcommand{\hds}{%
	\ensuremath{\alpha_{\mathrm{s}}}
}
\newcommand{\hdl}{%
	\ensuremath{\alpha_{\mathrm{l}}}
}

\newcommand{\dn}{%
	\ensuremath{\rho_{\gen}}
}
\newcommand{\ds}{%
	\ensuremath{\rho_{\mathrm{s}}}
}
\newcommand{\dmelt}{%
	\ensuremath{\rho_{\mathrm{m}}}
}
\newcommand{\dl}{%
	\ensuremath{\rho_{\mathrm{l}}}
}

\newcommand{\stc}{%
	\ensuremath{c_{\mathrm{p}}}
}
\newcommand{\stcs}{%
	\ensuremath{c_{\mathrm{p,s}}}
}
\newcommand{\stcl}{%
	\ensuremath{c_{\mathrm{p,l}}}
}

\newcommand{\bfd}{%
	\ensuremath{\delta_{\gen}}
}
\newcommand{\bfds}{%
	\ensuremath{\delta_{\mathrm{s}}}
}
\newcommand{\bfdl}{%
	\ensuremath{\delta_{\mathrm{l}}}
}
\newcommand{\dom}{%
	\ensuremath{\Omega_{\gen}}
}
\newcommand{\doms}{%
	\ensuremath{\Omega_{\mathrm{s}}}
}
\newcommand{\doml}{%
	\ensuremath{\Omega_{\mathrm{l}}}
}
\newcommand{\domr}{%
	\ensuremath{\Omega_{\mathrm{r}}}
}
\newcommand{\domd}{\Omega_{\mathrm{d}}}
\newcommand{\domdt}{\tilde\Omega_{\mathrm{d}}}

\newcommand{\inp}{%
	\ensuremath{u_{\gen}(t)}
}
\newcommand{\inpr}{%
	\ensuremath{\bs{u}(t)}
}
\newcommand{\inps}{%
	\ensuremath{u_{\mathrm{s}}(t)}
}
\newcommand{\inpl}{%
	\ensuremath{u_{\mathrm{l}}(t)}
}
\newcommand{\inpsr}{%
	\ensuremath{u_{\mathrm{s,r}}(t)}
}
\newcommand{\inplr}{%
	\ensuremath{u_{\mathrm{l,r}}(t)}
}
\newcommand{\inpcs}{%
	\ensuremath{u_{\mathrm{c,s}}(t)}
}
\newcommand{\inpcl}{%
	\ensuremath{u_{\mathrm{c,l}}(t)}
}
\newcommand{\uo}{%
	\ensuremath{\bs{u}(t)}
}
\newcommand{\ue}{%
	\ensuremath{\hat{\bs{u}}(t)}
}
\newcommand{\ub}{%
	\ensuremath{\bar{\bs{u}}(t)}
}
\newcommand{\ur}{%
	\ensuremath{\bs{u}_{\mathrm{r}}(t)}
}
\newcommand{\uf}{%
	\ensuremath{\bs{u}_{\mathrm{f}}(t)}
}
\newcommand{\uc}{%
	\ensuremath{\bs{u}_{\mathrm{c}}(t)}
}
\newcommand{\ut}{%
	\ensuremath{\tilde{\bs{u}}(t)}
}
\newcommand{\ubt}{%
	\ensuremath{\tilde{\bar{\bs{u}}}(t)}
}
\newcommand{\ud}{%
	\ensuremath{{\bs{\mu}}(t)}
}


\newcommand{\tderiv}[2]{%
	\ifthenelse{\isempty{#2}}{%
		#1(t)
	}{%
		{#1}^{(#2)}(t)
	}
}

\newcommand{\y}[1][]{%
	\ensuremath{%
		\tderiv{\bs{y}}{#1}
	}
}
\newcommand{\yr}[1][]{%
	\ensuremath{%
		\tderiv{\bs{y}_{\mathrm{r}}}{#1}
	}
}

\newcommand{\sysouts}[1][]{%
	\eta_{\gen}(t)
}
\newcommand{\sysout}[1][]{%
	\ifthenelse{\isempty{#1}}{%
		\bs{\eta}(t)
	}{%
		\bs{\eta}_{\mathrm{#1}}(t)
	}
}
\newcommand{\outmap}[2]{%
	\bs{h}
	\ifthenelse{\isempty{#1}}{%
	}{%
		_{\mathrm{#1}}
	}\left(#2\right)
}

\newcommand{\yn}{%
	\ensuremath{\bs{\eta}(t)}
}
\newcommand{\ye}{%
	\ensuremath{\hat{\bs{\eta}}(t)}
}
\newcommand{\yb}{%
	\ensuremath{\bar{\bs{\eta}}(t)}
}
\newcommand{\yt}{%
	\ensuremath{\tilde{\bs{\eta}}(t)}
}
\newcommand{\ynt}{%
	\ensuremath{\tilde{\bs{\eta}}(t)}
}
\newcommand{\yte}{%
	\ensuremath{\hat{\tilde{\bs{\eta}}}(t)}
}
\newcommand{\ybt}{%
	\ensuremath{\tilde{\bar{\bs{\eta}}}(t)}
}
\newcommand{\yd}{%
	\ensuremath{{\bs{\nu}}(t)}
}

\newcommand{\bs}[1]{\boldsymbol{#1}}
\newcommand{\partiell}[3][]{\frac{\partial^{#1}#2}{\partial{#3}^{#1}}}
\newcommand{\inner}[2]{\langle #1 \vert #2 \rangle}


\newcommand{\varphigz}[1]{%
	\ensuremath{\varphi_{\gen,#1}^N(z)}
}
\newcommand{\varphibgz}[1]{%
	\ensuremath{\bar{\varphi}_{#1}^N(\bar{z})}
}

\newcommand{\cgt}[1]{%
	\ensuremath{c_{\gen,#1}(t)}
}
\newcommand{\w}{%
	\ensuremath{w}
}
\newcommand{\wb}{%
	\ensuremath{\bar{w}^N}
}
\newcommand{\wi}{%
	\ensuremath{\w^N_i(t)}
}
\newcommand{\wbi}{%
	\ensuremath{\wb_{\gen, i}(t)}
}
\newcommand{\wbv}{%
	\ensuremath{\bar{\bs{\w}}^N_{\gen}(t)}
}
\newcommand{\wbvt}{%
	\ensuremath{\dot{\bar{\bs{\w}}}^N_{\gen}(t)}
}
\newcommand{\wbvs}{%
	\ensuremath{\bar{\bs{\w}}^N_{\mathrm{s}}(t)}
}
\newcommand{\wbvl}{%
	\ensuremath{\bar{\bs{\w}}^N_{\mathrm{l}}(t)}
}

\newcommand{\x}{%
	\ensuremath{\bs{x}(t)}
}
\newcommand{\xd}{%
	\ensuremath{\dot{\bs{x}}(t)}
}

\newcommand{\xn}{%
	\ensuremath{\bs{x}(t)}
}
\newcommand{\xnd}{%
	\ensuremath{\dot{\bs{x}}(t)}
}

\newcommand{\xe}{%
	\ensuremath{\hat{\bs{x}}(t)}
}
\newcommand{\xed}{%
	\ensuremath{\dot{\hat{\bs{x}}}(t)}
}

\newcommand{\xr}{%
	\ensuremath{\bs{x}_{\mathrm{r}}(t)}
}
\newcommand{\xrd}{%
	\ensuremath{\dot{\bs{x}_{\mathrm{r}}}(t)}
}

\newcommand{\xt}{%
	\ensuremath{\tilde{\bs{x}}(t)}
}
\newcommand{\xtd}{%
	\ensuremath{\dot{\tilde{\bs{x}}}(t)}
}

\newcommand{\xnt}{%
	\ensuremath{\tilde{\bs{x}}(t)}
}
\newcommand{\xntd}{%
	\ensuremath{\dot{\tilde{\bs{x}}}(t)}
}

\newcommand{\xte}{%
	\ensuremath{\hat{\tilde{\bs{x}}}(t)}
}
\newcommand{\xted}{%
	\ensuremath{\dot{\hat{\tilde{\bs{x}}}}(t)}
}

\newcommand{\xb}{%
	\ensuremath{\bar{\bs{x}}(t)}
}
\newcommand{\xbd}{%
	\ensuremath{\dot{\bar{\bs{x}}}(t)}
}

\newcommand{\xbt}{%
	\ensuremath{\tilde{\bar{\bs{x}}}(t)}
}
\newcommand{\xbtd}{%
	\ensuremath{\dot{\tilde{\bar{\bs{x}}}}(t)}
}

\newcommand{\Lobs}{%
	\ensuremath{\bs{L}(t)}
}
\newcommand{\cov}{%
	\ensuremath{\bs{\Pi}}
}
\newcommand{\covt}{%
	\ensuremath{\cov(t)}
}

\newcommand{\Vo}{%
	\ensuremath{V(\stateerr)}
}
\newcommand{\Vot}{%
	\ensuremath{\dot V(\stateerr)}
}
\newcommand{\Vg}{%
	\ensuremath{V_{\gen}(\stateerr)}
}
\newcommand{\Vgt}{%
	\ensuremath{\dot V_{\gen}(\stateerr)}
}
\newcommand{\Vs}{%
	\ensuremath{V_{\mathrm{s}}(\stateerr)}
}
\newcommand{\Vst}{%
	\ensuremath{\dot V_{\mathrm{s}}(\stateerr)}
}
\newcommand{\Vl}{%
	\ensuremath{V_{\mathrm{l}}(\stateerr)}
}
\newcommand{\Vlt}{%
	\ensuremath{\dot V_{\mathrm{l}}(\stateerr)}
}
\newcommand{\targetset}{%
	\mathcal{A}
}

\newcommand{\stateps}[1][N]{%
	\bs{\zeta}^{#1}_{\gen}(t)
}
\newcommand{\statepst}[1][N]{%
	\dot{\bs{\zeta}}^{#1}_{\gen}(t)
}

\newcommand{\statept}[1][N]{%
	\bar{\bs{\zeta}}^{#1}(t)
}
\newcommand{\stateptt}[1][N]{%
	\dot{\bar{\bs{\zeta}}}^{#1}(t)
}

\newcommand{\statef}[1][N]{%
	\bs{\chi}^{#1}(t)
}
\newcommand{\stateft}[1][N]{%
	\dot{\bs{\chi}}^{#1}(t)
}
\newcommand{\statefc}[1][N]{%
	\chi^{#1}_n(t)
}
\newcommand{\statefct}[1][N]{%
	\dot{\chi}^{#1}_n(t)
}
\newcommand{\statefs}[1][N]{%
	\tilde{\bs{\chi}}^{#1}(t)
}

\newcommand{\statemap}[2][N]{%
	\bs{\psi}^{#1}\left(#2\right)
}
\newcommand{\statemapi}[2][N]{%
	\bar{\bs{\psi}}^{#1}\left(#2\right)
}
\newcommand{\scalemat}[1][N]{%
	\bs{\mathcal{T}}^{#1}
}

\newcommand{\stateerr}{%
	\bs{\xi}
}
\newcommand{\stateerro}{%
	\bs{\xi}(z,0)
}
\newcommand{\stateerrz}{%
	\bs{\xi}(z,t)
}
\newcommand{\setdist}[1]{%
	\left\vert #1 \right\vert_{\targetset}
}

\newcommand{\rev}[1]{\textcolor{blue}{#1}}
\newcommand{\que}[1]{\textcolor{red}{#1}}

\section{Introduction}

The \gls{vgf} process is the most important technology for the production of bulk compound semiconductor
crystals like \gls{gaas} or \gls{inp} \cite{JURISCH2005283} which are especially used for manufacturing
high-power and high-frequency electronics as well as infrared light-emitting and laser diodes. For these
purposes the crystals have to meet high requirements with respect to their purity and structural
perfection. 

\begin{figure}
	\centering
	\def\svgwidth{.7\linewidth}
	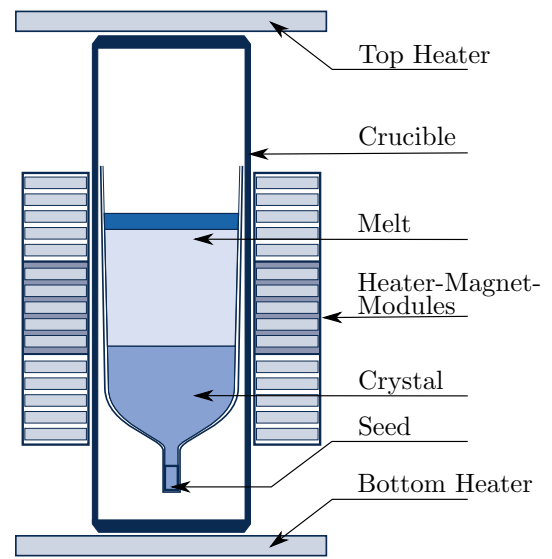
	\caption{Sketch of a \gls{vgf} crystal growth furnace.}
	\label{fig:PlantSketch}
\end{figure}%

The basic \gls{vgf} setup is shown in Figure \ref{fig:PlantSketch}: A crucible, usually made of
boron nitride and holding the material which is to be molten and then solidified, is surrounded by several
heaters. The whole setup is enclosed by a thick insulation. At the bottom of the crucible a seed crystal is
placed which defines the orientation of the crystal to be grown. After melting up the material in the crucible
(without destroying the seed crystal) the temperature field has to be adjusted and tracked by means of the
heat input of the heaters in such a way that the melt solidifies from the bottom to the top in a desired
manner. This means that 
a) the solid-liquid-interface (phase boundary) maintains a plain shape
and 
b) the growth rate (i.e., the velocity of the 
solid-liquid-interface) and the temperature gradient at the phase boundary are kept on a certain level throughout the whole process as they have been identified as crucial factors regarding the quality
of the grown crystal (see eg.~\cite{Vanhellemont2013} for an analysis regarding the related Czochralski process).
This solidification by using a travelling vertical temperature gradient is where the name of the 
process originates from.
To affect the process, a top and a bottom heater in form of plane disks, as well as three jacket heaters
in form of coils are mounted in the plant.

Due to its importance the improvement of this process is in focus of the scientific community,
resulting e.g.~in the application of external \glspl{tmf} \cite{FR2014,Dropka2013,DROPKA2014146,Dropka2014b}.
However, a topic that has not received much attention is the process control of this growth technique.  
This lack of coverage has two main reasons: 
Firstly, in-situ measurements from the growing crystal (e.g.~the shape
and position of the phase boundary or the growth rate) as a prerequisite for feedback control are
not available or not applicable in an industrial environment \cite{Wellmann2015}. 
Secondly, the coupled free boundary problems for crystal and melt form a so called two-phase
\gls{sp} \cite{Crank84} which is of nonlinear nature.

As is well known, the first issue can be tackled by an appropriate observer design which
has already been presented in \cite{Koga2016b} for the one-phase \gls{sp}.
However, regarding the implementation a simulation model is needed for the observer.
Since the simulation of solidification processes and therefore of \glspl{fbp} has been under 
investigation over the last decades, there are a lot of different numerical schemes like the 
Enthalpy\cite{Brusche2006}, Level-Set\cite{Chen1997} or Moving-Grid\cite{Beckett2001} method
available, to name just a few. However, being numerical schemes, identifying their variables
with a state-space representation for subsequent observer design is not straight forward. 

The second issue is broadly discussed in the framework of \glspl{dps}.
Making the assumption, that the temperature distribution in one phase is constant (which is often
justified due to its dominant spatial extent) yields the so called one-phase \gls{sp}.
Regarding this special case, results are lately available for the feedforward design 
\cite{DPRM03} using flatness-, as well as for feedback designs using enthalpy-
\cite{Petrus2012, Petrus2014}, geometrically- \cite{Maidi2014} or backstepping- \cite{Koga2016a} based
approaches.
Regarding the full problem, \cite{RWW04esta} extends the flatness-based motion planning to the
two-phase \gls{sp}, while \cite{Hinze2009} addresses the problem from the side of optimal control.
Concerning feedback, a direct extension of the approaches for the one-phase case is not feasible since for the two-phase case,
the coupling between the two \glspl{fbp} has to be taken into account.
In this context it is noteworthy that \cite{petrus2010} already states a Lyapunov-based control 
law for the two-phase \gls{sp} with actuation at one boundary.
However, according to our current knowledge there are no results available for the tracking control
of the two-phase \gls{sp} via output feedback concerning multiple inputs.

\subsection{Objective and structure of the paper} 
The main goal of this contribution is to introduce methods for tracking control of a one 
dimensional, two-phase \gls{sp} via output feedback (resp.~observer based state feedback) as a starting point for an improvement of
process control in the \gls{vgf} growth process. 

To reach this goal, the paper is structured as follows: 
In \autoref{sec:modelling} the distributed parameter model of the process is introduced.
\autoref{sec:feedforward} outlines a feedforward control design which is based on the flat
parametrisation of the solution by means of power series.
This feedforward subsequently serves as the source of a reference temperature profile.
Based on this, \autoref{sec:coll_feedback} introduces a collocated controller that tracks this
reference by utilizing state feedback.
Looking at the problem from another point of view and further exploiting the flatness property,
\autoref{sec:distr_feedback} presents a flatness based state feedback control.
This approach relies on a finite dimensional approximation of the system dynamics which is obtained
from the parametrisation in \autoref{sec:feedforward}.
To comply with the specific demands of the process, different variants of both control concepts are
introduced.
Since all designs depend on state measurements to some extend, in \autoref{sec:observer} a lumped
observer for the flat system approximation is shown.
\autoref{sec:results} presents
simulation results for the different control setups using state and output feedback.
Finally, a summary and an outlook to further work is given.

\section{Modelling of the VGF process}
\label{sec:modelling}

\begin{figure}
	\centering
	\def\svgwidth{\linewidth}
	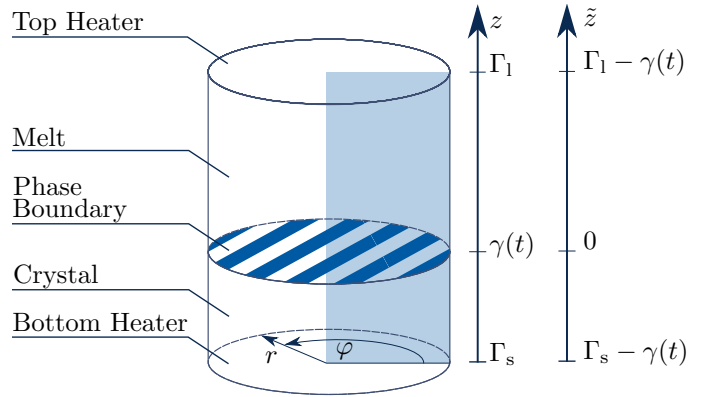
	\caption{Schematics of the cylindrical coordinate system $(r, \varphi, z, t)$, a meridional plane (blue)
		and the shifted coordinate $\zt = z - \zi$.}
	\label{fig:system}
\end{figure}%

In this section a one dimensional distributed parameter model
of the \gls{vgf} process plant 
is derived. 
For this purpose, the following simplifications are made:
The crucible geometry is approximated by a cylinder.
Thus, the distribution of the system temperature $T$ in the crucible depends on the time $t$ and the
cylindrical coordinates, given by radius $r$, angle $\varphi$ and height $z$, as depicted in
\autoref{fig:system}.
Furthermore, any convective effects in the melt are neglected.
This is reasonable due to the dominating heat transport by diffusion. 
Beyond, making use of the fact that the plant itself is rotationally symmetric to the longitudinal axis,
the model can be reduced to a meridional plane of the crucible by taking the average over the
angular coordinate $\varphi$.
In addition, the lateral heaters are assumed to be used as active isolation, avoiding any heat loss in
radial direction and therefore rendering the phase boundary a horizontal line.
Hence, averaging over the radius $r$ allows further reduction of the domain to a line whose
boundaries represent the bottom and top of the crucible at $z=\zbs$ and $z=\zbl$, respectively.
Summarising, the system temperature is given by $T(z,t)$ for $\zbs \le z \le \zbl$ and $t>0$ while
the phase boundary is given by $\zi \in (\zbs, \zbl)$.

This leads to the one dimensional nonlinear heat equation \cite{Rota84}
\begin{align}
	\begin{split}
		\partiell{}{t}\Big(\rho(\tnz)\stc(\tnz) \tnz \Big)=\\
		\qquad \partiell{}{z}\left(k(\tnz) \partiell{}{z} \tnz \right)
		,z \in (\zbs, \zbl) \setminus \left\{\zi\right\}
	\end{split}
	\label{eq:nonlin_heat_eq}
\end{align}
with the density $\rho$, the specific heat capacity \stc and $k$ the thermal conductivity
being temperature-dependent.
Assuming piecewise
constant parameters for the solid and the liquid phase it is possible to decompose the nonlinear system
\eqref{eq:nonlin_heat_eq} into two \glspl{fbp} for the temperatures $\ts(z,t)$ and $\tl(z,t)$:
\begin{subequations}
	\label{eq:lin_heat_eqs}
	\begin{align}
		\partial_t \ts(z, t) &= \hds\partial_z^2 \ts(z,t),\;z \in \doms = (\zbs, \zi)\\
		\label{eq:lin_heat_eq_s_bc1}
		\hcs\partial_z \ts (\zbs, t) &= \bfds\inps\\
		\ts(\zi, t) &= \tm \\
		\partial_t \tl(z, t) &= \hdl\partial_z^2 \tl(z,t),\; z \in \doml = (\zi, \zbl)\\
		\label{eq:lin_heat_eq_l_bc1}
		\hcl\partial_z \tl (\zbl, t) &= \bfdl\inpl\\
		\tl(\zi, t) &= \tm.
	\end{align}
\end{subequations} 
Herein, the index ``s'' denotes the solid and the index ``l'' the liquid phase.
The heat flows $\inps$ and $\inpl$ at the bottom and the top boundary are considered as system inputs 
with the orientation factors $\bfds=-1$ and $\bfdl=1$. 
The partial derivative of the quantity $\torigz$ with respect to $z$ or $t$ is denoted by
$\partial_z \tnz$ or $\partial_t \tnz$.
Finally, $\hds = \frac{\hcs}{\ds\stcs}$ and $\hdl = \frac{\hcl}{\dl\stcl}$ denote the thermal 
diffusivities.

Due to the moving phase boundary latent heat is released by the solidification process. 
This effect can be modelled by the \acrlong{sc} \cite{stefan_cond}
\begin{equation}
	\dmelt L\vi = \hcs\partial_z\ts(\zi, t) - \hcl\partial_z\tl(\zi, t)
	\label{eq:stefan_cond}
\end{equation}
with the density of the melt at melting temperature $\dmelt$ and the specific latent heat $L$.

Together, the equations \eqref{eq:lin_heat_eqs} and \eqref{eq:stefan_cond} form the two-phase 
\gls{sp} whose state is given by
\begin{equation}
	\bs{x}(\cdot, t) = \begin{pmatrix}
		\torig(\cdot, t) \\ \zi
	\end{pmatrix} \in X = 
		L_2(\Omega)
		\times (\zbs, \zbl)
	\label{eq:state}
\end{equation}
with $\Omega = [\zbs, \zbl]$.
Note that the PDE-ODE system defined by 
\eqref{eq:lin_heat_eqs} and \eqref{eq:stefan_cond} is inherently nonlinear since the domains of
\eqref{eq:lin_heat_eq_s_bc1} and \eqref{eq:lin_heat_eq_l_bc1} depend on the state variable $\zi$.
Furthermore, with the system boundaries admitting access for measurements, the system output
is given by
\begin{equation}
	\sysout = \outmap{}{\x} = \begin{pmatrix}
		\torig(\zbs, t)\\
		\torig(\zbl, t)
	\end{pmatrix}.
	\label{eq:output}
\end{equation}

Exploiting the identical structure of the diffusion equations, the following sections will 
-- where applicable --  resort to discuss merely one generic temperature distribution 
$\torigg(z,t)$ for $z \in \dom$
with $\circ$ to be replaced by the indices $s$ or $l$ depending on the considered domain.
%


\section{Feedforward control}
\label{sec:feedforward}

This section gives a short recap of a feedforward control design which was presented in \cite{DPRM03} for
the one-phase and in \cite{RWW03md} for the two-phase case to which the reader is kindly directed
for further details.

To eliminate the temporal dependency in the boundary conditions of \eqref{eq:lin_heat_eqs}, the 
coordinate transform
\begin{equation}
	\ttransgz = \toriggz \quad \text{with} \quad \zt \coloneqq z - \zi
	\label{eq:coord_trans}
\end{equation}
is introduced.
As a consequence the phase boundary is shifted into the origin of a new, moving reference frame as it
can be seen on the right-hand side of \autoref{fig:system}.
The resulting\footnote{For details see \ref{app:trafo_mrs}.} system is given by
\begin{subequations}
	\begin{align}
		\label{eq:transformed_pde}
		\partial_t \ttransgz &= \hd\partial_{\zt}^2\ttransgz + \vi \partial_{\zt} \ttransgz \\
		\label{eq:transformed_bc_0}
		\hc\partial_{\zt} \ttrans(\zbt, t) &= \bfd\inp \\
		\label{eq:transformed_bc_1}
		\ttransg(0, t) &= \tm \\
		\label{eq:transformed_sc}
		\vi &=
		\frac{1}{L \dmelt}\left(
			\hcs\partial_{\zt}\ttranss(0, t) -\hcl\partial_{\zt}\ttransl(0, t)
		\right)
	\end{align}%
	\label{eq:transformed_system}%
\end{subequations}
where $\zbt \coloneqq \zb - \zi$.

By expressing the solution $\ttransgz$ of \eqref{eq:transformed_system} in terms of a power series
in $\zt$:
\begin{equation}
	\ttransgz= \sum\limits_{i=0}^\infty \cgt{i}\frac{\zt^i}{i!} \:,
	\label{eq:power_series}
\end{equation}
plugging it into \eqref{eq:transformed_pde} and comparing the coefficients of like powers of
$\zt$ the recursion formula
\begin{equation}
	\cgt{i+2} = \frac{1}{\hd}\big(
		\partial_t \cgt{i}
		- \vi \cgt{i+1}
	\big)
	\quad i=0,\dotsc,\infty
	\label{eq:coeff_iteration}
\end{equation}
is obtained.
A closer examination of \eqref{eq:power_series} shows that the following holds for the initial
coefficients
:
\begin{alignat}{2}%
	\cgt{0} &= \ttransg(0, t) = \tm ,
	&\quad \cgt{1} &= \partial_{\zt} \ttransg(0, t) \:.
	\label{eq:c1}
\end{alignat}%
By utilizing the \acrlong{sc} \eqref{eq:transformed_sc} solved for the gradient at the liquid side
\begin{equation}
\partial_{\zt}\ttransl(0, t) = \frac{1}{\hcl}\left(\hcs\partial_{\zt}\ttranss(0, t) 
	- \dmelt L\vi\right)
\end{equation}
it follows that the solution for both phases can be expressed by the gradient in the solid
$\partial_{\zt}\ttranss(0, t)$ and the growth rate $\vi$.
Thus, the system \eqref{eq:lin_heat_eqs} is differentially flat with a flat output
\begin{equation}
	\bs{y}(t) 
	= \begin{pmatrix} y_1(t) \\ y_2(t)\end{pmatrix}
	= \begin{pmatrix}
			\partial_{\zt}\ttranss(0, t) \\
			\zi
		\end{pmatrix}.
	\label{eq:flat_output}
\end{equation}
Assuming convergence and truncating the series \eqref{eq:power_series} at an order $N$, the mapping
\begin{align}
	\begin{split}
		\tilde{\bs{T}}^N(\tilde z, t) = \begin{pmatrix}
			\ttranss^N(\tilde z, t) \\ \ttransl^N(\tilde z, t)
		\end{pmatrix} = \bs{\Theta}^N\big(
			&\partial_{\zt}\ttranss(0, t),
			\dotsc,
			\partial_{\zt}\ttranss^{(\alpha_1)}(0, t),\\
			&\zi, \dotsc, \gamma^{(\alpha_2)}(t)\big)
	\end{split}
	\label{eq:flat_state_map}
\end{align}
where $\alpha_1 = \floor{N/2}-1$ and $\alpha_2 = \floor{N/2}$ can be formulated.
%
Finally, choosing the trajectories for the components of $\bs{y}(t)$ as
\begin{equation}
	y(t) = y_A +(y_B- y_A)\phi\left(\frac{t}{\vartheta}\right)
	\label{eq:flat_traj}
\end{equation}
where 
\begin{equation}
	\phi(\tau) = \frac{1}{2}\left(1 + \tanh\left(\frac{2(2\tau - 1)}
		{(4\tau(1-\tau))^\sigma}\right)\right)
	\label{eq:smooth_func}
\end{equation}
is of Gevrey-order $\alpha_{\textup{G}} = 1+\frac{1}{\sigma}$ \cite{Gevrey1918} convergence of the
scheme can be shown for $\alpha_{\textup{G}}\le2$ \cite{RWW04esta}.

Hence, by utilizing the boundary condition \eqref{eq:transformed_bc_0} the inputs
$\inps$
and
$\inpl$
are obtained by evaluating the temperature profile which is given the mapping 
\eqref{eq:flat_state_map}.
This yields the input map
\begin{equation}
	\inpr = \bs{\Phi}^N(y_1(t), \dotsc, y_1^{(\alpha_1)}(t), 
		y_2(t), \dotsc, y_2^{(\alpha_2)}(t))
	\label{eq:flat_input_map}
\end{equation}
given in terms of the flat output.

Summarising, by prescribing reference trajectories $\yr$ for the flat output $\y$, a reference
temperature distribution $\trefz$ and input $\ur$ can be computed.
These results will serve as the basis for the control design in the next section.


\section{Collocated feedback}
\label{sec:coll_feedback}

This section will present an approach for the collocated tracking control of the
two-phase \gls{sp}.
To do this, different control errors are compared concerning their suitability
concerning the \gls{vgf} process and the stability of the closed loop is
investigated.
In detail, the profile $\trefz$, introduced in the previous section, will be utilised as control
reference while the input trajectories $ \ur$ are not needed.
Note that all considerations in this section are conducted in the original spatially fixed coordinates.

\subsection{Error definitions}
\label{ssec:error_system}
Choosing the distributed temperature error
\begin{equation}
	\terrn(z, t) \coloneqq \torigz - \trefz
	\label{eq:fixed_error}
\end{equation}
as deviation of the system temperature $\torigz$ from the reference profile $\trefz$ over the
complete spatial domain as in \cite{petrus2010} seems natural at a first glance.
However, this is not useful in order to meat the technological requirements:
Due to the biphasic character of the system, convergence of $\terrn(z,t)$ into the origin implies
convergence of the phase boundary position error
\begin{equation}
	\dzi = \zi - \zir \:.
	\label{eq:zi_error}
\end{equation}
However, the  quality of the crystal does not depend on the position of the phase boundary but rather on its velocity.
Moreover, it seems more natural to compare the temperature of like phases only.
As a consequence the error is defined on the basis of a shifted reference trajectory
\begin{equation}
	\terrz \coloneqq \torigz - \treft
	\label{eq:shifted_error}
\end{equation}
in combination with \eqref{eq:zi_error}.
To do so, however, the planning for the reference profile $\trefz$ has to be carried out on an extended
spatial domain. More precisely, for any admissible $\dzi$, the profile $\treft$ must not be
evaluated outside of its domain. Considering the plant properties, it follows that
$\dzi \in (\zbs-\zbl, \zbl-\zbs)$. Therefore, a feasible domain for $\trefz$ would be
$(z, t) \in \domr\times\mathbb{R}_+$, with $\domr = [2\zbs-\zbl, 2\zbl-\zbs]$.
Analysing\footnote{Details in \ref{app:trafo_err}.} \eqref{eq:shifted_error}, 
the phase dependent tracking error $\terrgz = \toriggz - \trefgt$ is governed by
\begin{subequations}
	\begin{align}
		\label{eq:error_system_pde}
		\partial_t \terrgz &= \hd\partial_z^2 \terrgz +  \dvi \partial_z\trefgt
		\nonumber\\&\quad\quad z \in \dom\\
		\label{eq:error_system_bc1}
		\partial_z \terrg(\zb, t) &= \frac{\bfd}{\hc} u(t)
		- \partial_z \trefg(\zb - \dzi, t) \\
		\label{eq:error_system_bc2}
		\terrg(\zi, t) &= 0 \\
		\label{eq:error_system_sc}
		\dvi &= \frac{1}{L \dmelt}\left(\hcs\terrs(\zi, t) - \hcl\terrl(\zi, t)\right)\:.
	\end{align}%
	\label{eq:error_system}%
\end{subequations}%
Herein, the second term in the rigth-hand side of \eqref{eq:error_system_bc1} can be understood as
a feedforward part.
However, it does not coincide with the reference input $\ur$ due to the spatial shift.
Finally, the error state is given by
\begin{equation}
	\stateerrz = \begin{pmatrix} \terrgz  \\ \dzi \end{pmatrix} \in X .
	\label{eq:err_state}
\end{equation}

\subsection{Control law}
A very intuitive way to manipulate the system is to convert \eqref{eq:error_system_bc1} into
%
\begin{equation}
	\partial_z\terrg(\zb, t) = -\bfd\cgaing\terrg(\zb, t) \:.
	\label{eq:auto_bc_1}
\end{equation}
This leads to the feedback law 
\begin{equation}
	\inp = \frac{\hc}{\bfd} \partial_z\trefgt - \cgaing \hc \terrg(\zb, t).
	\label{eq:bc_feedback}
\end{equation}
Albeit reasonable in its composition, the control law only honours the boundary error at $z = \zb$.
Hence, further analysis is required to ensure convergence of $\dvi$.

\subsection{Stability analysis}
Although the framework which will be applied here was already laid out in \cite{Zubov1964} for
finite dimensional systems, the nomenclature, used in the following is borrowed from 
\cite{Mironchenko17} due to its application for the infinite dimensional case.
Keeping in mind that the tracking of the growth velocity $\vi$ is more important than the exact
adjustment of the boundary position $\zi$, it is apparent that to obtain the desired results,
$\stateerrz$ may not necessarily converge into the origin but rather into a compact subset of the
state space, given by
\begin{equation}
	\targetset \coloneqq \left\{(\stateerr, \gamma)^T \in X \vert\,
	\stateerr = 0 \right\}
	\label{eq:target_set}
\end{equation}
with $X$ from \eqref{eq:state}.
Assuming that $\stateerrz \in \targetset$ for some $t$, \eqref{eq:error_system_sc} yields that $\dvi = 0$.
Using this information, \eqref{eq:error_system_pde} gives $\partial_t \terrz = 0$.
Thus $\stateerro \in \targetset \implies \stateerrz \in \targetset \quad\forall t \ge 0$,
rendering $\targetset$ an invariant set of the system \eqref{eq:error_system}.
Furthermore let the distance of an element $\bs{x} \in X$ to $\targetset$ be given by 
$\setdist{\bs{x}} \coloneqq \min\left\{\norm{\bs{x}-\bs{y}}_{X} \vert\, \bs{y} \in \targetset\right\}$
and consider the function classes:
\begin{align*}
	\mathcal{K} &\coloneqq \{f:\mathbb{R}_+ \mapsto \mathbb{R}_+ 
		\vert\, f(0) = 0, f \text{ is continuous} \\ &\qquad \text{ and strictly increasing}\} \\
		\mathcal{K_{\infty}} &\coloneqq \left\{f \in \mathcal{K}\vert\, f \text{ is unbounded}\right\} \:.
\end{align*}
As stated in \cite{Mironchenko17}, if there exists a Lyapunov function $\Vo$\footnote{
	For better readability, in the following the arguments of $\stateerrz$ are omitted.
}, so that
\begin{subequations}
	\begin{align}
		\label{eq:pos_def}%
		a_1(\setdist{\stateerr}) \le V(\stateerr) &\le a_2(\setdist{\stateerr})\quad\forall \stateerr \in X \\
		\label{eq:neg_def}%
		\dot V(\stateerr) &\le -b(\setdist{\stateerr}) \quad\forall \stateerr \in X
	\end{align}%
	\label{eq:ugas_cond}%
\end{subequations}%
with $a_1, a_2 \in \mathcal{K_{\infty}}$ and $b \in \mathcal{K}$ hold, 
the system \eqref{eq:error_system} is uniformly globally asymptotically stable with respect 
to $\targetset$. 
For this purpose, the Lyapunov function candidate
\begin{equation}
	\Vo = \frac{1}{2} \int\limits_{\zbs}^{\zbl} \terrzsq \, \textup{d}z
	\label{eq:ljap_cand}
\end{equation}
is used, which fulfils condition \eqref{eq:pos_def}.
Furthermore, in the first part of \ref{app:a} it is shown that for a simplified reference profile $\trefz$,
the candidate \eqref{eq:ljap_cand} satisfies \eqref{eq:neg_def}. 
Thus, rendering it a Lyapunov function for \eqref{eq:error_system} with respect to $\targetset$.
Softening those demands on $\trefz$ is possible but leads to stricter requirements for the phase
boundary error $\dzi$, which are again hard to show for the general case.
The detailed steps are given in in the second part of \ref{app:a}. 
However, simulation results show that the system state $\stateerr$ converges to $\targetset$ 
for non-trivial reference profiles, too.



\section{Distributed Feedback}
\label{sec:distr_feedback}
Exploiting the fact that parametrisation of the system \eqref{eq:lin_heat_eqs}-\eqref{eq:stefan_cond}
which is introduced in \autoref{sec:feedforward} is differentially flat \cite{DPRM03, RWW04esta} with the flat
output \eqref{eq:flat_output},
this property can be used to design a feedback in a straight-forward fashion without 
the explicit computation of a reference temperature profile.

\subsection{System State}
\label{ssec:flat_sys}
In \cite[Ch.~5]{MeurerDiss} a state space representation for a diffusion equation
is given by means of the series coefficients of a power approximation.
Therein, the components of a new state $\stateps \coloneqq (\cgt{1}, \dotsc, \cgt{N})^T$
belong to an approximation of an order $N$.
Instead of extending this approach by combining the coefficients of the solid and liquid
approximations into an extended state vector 
$\statept = \left(\bs{x}_{\mathrm{s}}^{N^T}(t), \bs{x}_{\mathrm{l}}^{N^T}(t) \right)^T$
one may directly use the appropriate derivatives of the flat output.
Hence, the state components $\statefc, 1\le n\le M$ with $M = (\alpha_1 + \alpha_2 +2)$
in flat coordinates constitute the state vector $\statef \in \mathbb{R}^{M}$ which reads:
\begin{equation}
	\statef = 
		\left(
		y_{1}(t),
		\dotsc, 
		y_{1}^{(\alpha_1)}(t),
		y_{2}(t),	
		\dotsc, 
		y_{2}^{(\alpha_2)}(t)
	\right)^T.
	\label{eq:flat_state}
\end{equation}
Herein, the required derivatives can be obtained from the iterated recursion formulas for both phases,
cf.~\eqref{eq:coeff_iteration}, for clarity condensed in the map
\begin{equation}
	\statef = \statemap{\statept}.
	\label{eq:flat_coeff_map}
\end{equation}
Thus, examining the components of the derivative $\stateft$ two integrator chains become apparent:
\begin{align}
	\statefct = \begin{cases}
		y_1^{(n)}(t)				& \text{for }1 \le n \le \alpha_1 + 1 \\
		y_2^{(n-(\alpha_1+1))}(t) 	& \text{for }\alpha_1 + 1 < n \le M\:.
	\end{cases}
	\label{eq:flat_state_dt}
\end{align}
Herein, the yet unknown derivatives $y_1^{(\alpha_1+1)}$ and $y_2^{(\alpha_2+1)}$ can be obtained
from an extended version of \eqref{eq:flat_coeff_map} by using the extended coefficient state 
$\statept[N+1]$:
\begin{equation}
	\statef[N+1] = \statemap[N+1]{\statept[N+1]}.
	\label{eq:flat_coeff_map_ext}
\end{equation}
However this mapping requires the coefficients $c_{\mathrm{s}, N+1}(t)$ and 
$c_{\mathrm{l}, N+1}(t)$.
Fortunately, these can be acquired from the respective boundary conditions of both phases,
cf.~\eqref{eq:transformed_bc_0}, after inserting the series expansion
\begin{equation}
	\cgt{N+1} = \frac{N!}{\zb^N}\hc\left(\bfd \inp 
	- \sum\limits_{i=0}^{N-1} \cgt{i+1}\frac{\zt^i}{i!}\right)
	\label{eq:coeef_np1}
\end{equation}
wherein the coefficients $\cgt{1}$ to $\cgt{N}$ can readily be computed from $\statef$.
According to \eqref{eq:output}, the outputs of each phase are given by
\begin{equation}
	\sysouts = \ttransg(\zbt, t) =  \sum\limits_{i=0}^{N} \cgt{i}\frac{\zbt^i}{i!}\:.
	\label{eq:meas_eq}
\end{equation}
Hence, by using $\statemapi{\cdot}$, the inverse of the map \eqref{eq:flat_coeff_map}, the output can be written as
\begin{equation}
	\sysout[flat] = \outmap{flat}{\statemapi{\statef}}.
	\label{eq:flat_out_eq}
\end{equation}

\subsection{Feedback Design}
Regarding $y_1(t)$ and $y_2(t)$ as the outputs of the system,
the tracking errors 
\begin{equation}
	\varepsilon_j(t) = y_j(t) - y_{j,\mathrm{r}}(t),\quad j=1,2
	\label{eq:tracking_error}
\end{equation}
are defined.
Hence, the decoupled linear error dynamics
\begin{subequations}
	\begin{align}
		\label{eq:error_dynamic_1}%
		\varepsilon_{1}^{(\alpha_1+1)}(t) &=
		- \sum\limits_{i=0}^{\alpha_1} \kappa_{1,i} \varepsilon_1^{(i)}(t) 
			\\
		\label{eq:error_dynamic_2}%
		\varepsilon_{2}^{(\alpha_2+1)}(t) &=
		- \sum\limits_{i=0}^{\alpha_2} \kappa_{2,i} \varepsilon_2^{(i)}(t) 
	\end{align}%
	\label{eq:error_dynamics}%
\end{subequations}%
are prescribed by choosing appropriate coefficients $\kappa_1$ and $\kappa_2$.
Defining the new inputs 
$v_1(t) \coloneqq y_1^{(\alpha_1+1)}(t)$ and $v_2(t) \coloneqq y_2^{(\alpha_2+1)}(t)$,
by using the inverse of \eqref{eq:flat_coeff_map_ext}
\begin{equation}
	\statept[N+1] = \statemapi[N+1]{\statef, (v_1(t), v_2(t))^T},
	\label{eq:ext_ceoff_flat_map}
\end{equation}
the extended series coefficient set can be computed.
Lastly, evaluation of \eqref{eq:flat_input_map} yields the control input $\uo$.

%
%
However, this design inherits the problem that an already grown crystal will be remelted if the 
measured interface position is ahead of the reference.
Therefore, in view of the shifted error system \eqref{eq:shifted_error}, the pair ($y_1(t), \dot y_2(t))$
may be regarded as the output of a modified system, yielding
$\tilde \varepsilon_2(t) = \dot y_2(t) - \dot y_{2,\mathrm{r}}(t)$ as well as the dynamics
\begin{equation}
	\tilde\varepsilon_{2}^{(\alpha_2)}(t) =
	- \sum\limits_{i=0}^{\alpha_2-1} \tilde\kappa_{2,i} \tilde\varepsilon_2^{(i)}(t) 
	\label{eq:error_dynamic_2_mod}%
\end{equation}
instead of \eqref{eq:error_dynamic_2}.
Thus, a modified virtual input can be stated as
$\tilde v_2(t) \coloneqq \tilde\varepsilon_{2}^{(\alpha_2)}(t) + y^{(\alpha_2+1)}_{2,\mathrm{r}}(t)$
which, by using \eqref{eq:ext_ceoff_flat_map} and \eqref{eq:flat_input_map} with $\tilde v_2(t)$
instead of $\tilde v_2(t)$, yields the modified control input $\tilde{\bs{u}}(t)$.

\section{Observer Design}
\label{sec:observer}

This section performs an observer design as shown in \cite{Meurer2005}, however in this case
based on the flat system state \eqref{eq:flat_state}.
The estimated system with the state $\xe$ is given by the copy
\begin{subequations}
	\begin{align}
		\xed &= \bs{f}(\xe, \uo) + \Lobs\yb \\
		\ye &= \bs{h}(\xe)
	\end{align}
	\label{eq:ext_luen}%
\end{subequations}
with $\Lobs$ to be chosen later on and $\yb = \ye - \yn$.
%
%
The plant model \eqref{eq:ext_luen} is extended in the following way:
\begin{equation}
	\xnd = \bs{f}(\xn, \uo + \ud),\quad \yn = \bs{h}(\xn + \yd) \:.
	\label{eq:nonlin_ss_dist}
\end{equation}
Herein, $\ud$ and $\yd$ represent disturbances acting on the system input and output, respectively.
%
Furthermore, denoting $\xb = \xe - \xn$,
the observer error dynamics
\begin{subequations}%
	\begin{align}%
		\begin{split}%
			\xbd &=  \bs{f}(\xe, \uo) 
			+\Lobs\yb
			\\&\quad
			- \bs{f}(\xn, \uo - \ud)
		\end{split}\\
		\yb &= \bs{h}(\xe) - \bs{h}(\xn - \yd)
	\end{align}%
	\label{eq:obs_err}%
\end{subequations}
is obtained.
In the following, the computation of $\Lobs$ will be performed on a linearisation of
\eqref{eq:obs_err} along the reference trajectory $\bs{y}_{\mathrm{r}}(t)$, given by:
\begin{subequations}
	\begin{align}
		\label{eq:lin_lumped_state}
		\xbd &= \bs{A}(t)\xb - \bs{B}(t)\ud +  \Lobs\yb\\
		\label{eq:lin_lumped_out}
		\yb &= \bs{C}(t)\left(\xb - \yd \right).
	\end{align}%
	\label{eq:lin_lumped}%
\end{subequations}%
%
By defining the cost functional
\begin{equation}
	J = \bar{\bs{x}}^T(0)\bs{S}\bar{\bs{x}}(0)
	+ \int\limits_0^t \bs{\mu}^T(t)\bs{R}\ud + \bar{\bs{\eta}}^T(t)\bs{Q}\yb \,\textup{d}t
	\label{eq:lqe_functional}
\end{equation}
where $\bs{S} \in \mathbb{R}^{M\times M}$ and $\bs{R}, \bs{Q} \in \mathbb{R}^{2\times 2}$ denote
penalties concerning the initial error as well as the disturbances on input and output, respectively.
As \cite[Th.~40, p.378]{sontag98} states, using the solution $\covt \in \mathbb{R}^{M\times M}$ of
the \gls{frde}
\begin{subequations}
	\begin{align}
		\begin{split}
			\dot{\cov}(t) &=
			\covt\bs{A}^T(t)
			+\bs{A}(t)\covt \\
			&\quad- \covt \bs{C}^T(t)\bs{Q}\bs{C}(t)\covt\\
			&\quad+ \bs{B}(t)\bs{R}^{-1}\bs{B}^T(t)
		\end{split}
		\label{eq:frde_ode}
		\intertext{with the initial condition}
		\cov(0) &= \bs{S}^{-1},
		\label{eq:frde_ic}%
	\end{align}%
	\label{eq:frde}%
\end{subequations}
the choice
\begin{equation}
	\Lobs \coloneqq -\covt\bs{C}^T(t)\bs{Q}
	\label{eq:def_l}
\end{equation}
yields the optimal estimation
for \eqref{eq:lin_lumped} regarding \eqref{eq:lqe_functional}.
%
%
Note that the solution of \eqref{eq:frde} can be done in advance.


\section{Results}
\label{sec:results}
The theoretical results of the previous sections will now be evaluated by simulations.
A \gls{fem} approximation using the boundary-immobilisation method
will serve as a simulation model to compare the different feedback designs on a process oriented benchmark
from the \gls{vgf} process.
The corresponding parameters are given in \autoref{tab:sim_params}.

\begin{table} 
	\centering
	\caption{Parameters of the System}
	\setlength{\tabcolsep}{0.3em}
	\begin{tabular}{@{}lcccl@{}} 
		\toprule
		Name	&Symbol	&Value (s/l)	&Unit \\
		\midrule
		Spec.~heat cap.	&$\stc$	&\num{423,59} /	\num{434}		&\si{\joule\per\kilogram\per\kelvin} \\
		Therm.~conduct.		&$k$		&\num{7,17} /	\num{17,8}		&\si{\watt\per\metre\per\kelvin} \\
		Therm.~diffus.		&$\hds$		&\num{3,27e-6}/					&\si{\square\metre\per\second} \\
		&$\hdl$	& \num{7,19e-6}	&\si{\square\metre\per\second} \\
		Densities					&$\ds$		&\num{5171,24} /	&\si{\kilogram\per\cubic\metre} \\
		&$\dl$	&\num{5702.37}	&\si{\kilogram\per\cubic\metre} \\
		&$\dmelt$	&\num{5713,07}	&\si{\kilogram\per\cubic\metre} \\
		Melting temp.		&$\tm$		&\num{1511,15}					&\si{\kelvin} \\
		Spec.~latent heat	&$L$		&\num{668,5e3}					&\si{\joule\per\kilogram} \\
		Left Boundary			&$\zbs$		&\num{0}						&\si{\metre}\\
		Right Boundary			&$\zbl$		&\num{0,4}						&\si{\metre}\\
		Feedf.~Appr.~Order		&$N_{\mathrm{ff}}$			&\num{10}						\\
		Obs.~Approx.~Order		&$N_{\mathrm{ob}}$		&\num{5}						\\
		Cont.~Approx.~Order		&$N_{\mathrm{fb}}$		&\num{5}						\\
		Inp.~Weight Mat.~			&$\bs{R}$		
		&\num{1e-4}$\bs{I}_2$
			&\si{\metre\tothe4\per\square\watt\per\second}\\
		Out.~Weight Mat.			&$\bs{Q}$		
			&\num{1e-4}$\bs{I}_2$
			&\si{\per\square\kelvin\per\second}\\
		Inp.~ dist.
			&\multicolumn{2}{c}{$\ud = \mathcal{N}\left(0, 100\right)$
				\tablefootnote{\phantomsection\label{fn:unit_note}The given unit refers to the 
						second argument}
			}
			&\si{\kilo\watt\squared\per\metre\tothe{4}}\\
		Output dist.~
			&\multicolumn{2}{c}{$\yd = \mathcal{N}\left(0, 10\right)$
				\footref{fn:unit_note}
			}
			&\si{\kelvin\squared}\\
		Init.~Weight Mat.			&$\bs{S}$
		&$\num{1e-3}\bs{I}_{5}$
		&$\si{\square\metre\per\square\kelvin}$\\
		Sim.~Disc.~Nodes		&$N_{\mathrm{fem}}$			&\num{41}					\\
		\bottomrule
	\end{tabular}
	\label{tab:sim_params}
\end{table}

\subsection{Setup and feedforward}
\label{ssec:results_ff}
\begin{figure}
	\centering
	\resizebox {\linewidth} {!} {
		\input{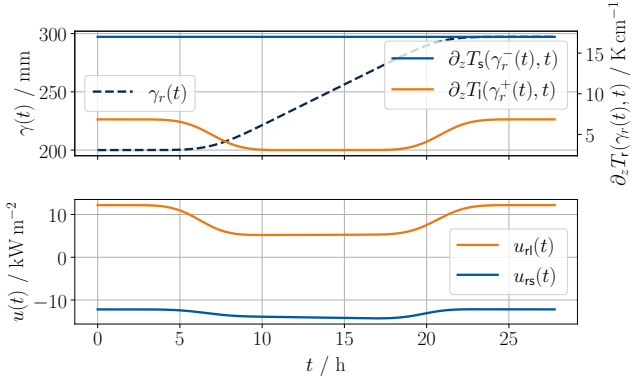}
	}
	\caption{Reference trajectories for gradients and phase boundary (top)
	as well as the generated heater trajectories for the system inputs (bottom)
	of the feedforward control.}
	\label{fig:trajectories}
\end{figure}

For the trajectory planning, the following initial situation is assumed: 
The phase boundary 
is resting ($\dot \gamma(0) = \SI{0}{\metre\per\second}$)
at $\gamma(0) = \SI{0,2}{\metre}$.
Furthermore, as a result of a previous step (as shown in \cite{RWW04esta}) a gradient of
$\partial_t\ts(\gamma(0), 0) = \SI{17}{\kelvin\per\centi\metre}$ has been established
at the solid side of the phase boundary.
Now, the growth process is performed by prescribing $\zir$.
\autoref{fig:trajectories} shows the generated trajectories for 
$\partial_t\trefs(\zir, t)$ and $\zir$ as well as the calculated system inputs 
$\inps$ and $\inpl$.

\begin{figure}
	\centering
	\resizebox {\linewidth} {!} {
		\input{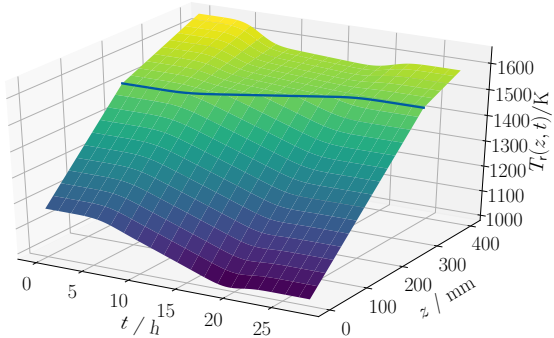}
	}
	\caption{Calculated reference temperature profile $\trefz$ with the reference
		phase boundary trajectory $\zir$ (blue).
	}
	\label{fig:ref_profile}
\end{figure}

\subsection{Feedback}
To emulate a real growth process, an initial error of $\zierro = \SI{100}{\milli\metre}$ and 
$\vierro = \SI{-3}{\milli\metre\per\hour}$ is introduced to the test-setup.
To gain an extensive overview, two versions of the collocated controller from \autoref{sec:coll_feedback}
are evaluated, one using the fixed error definition from \eqref{eq:fixed_error} and one using the
shifted error from \eqref{eq:shifted_error}.
Furthermore, two variants of the distributed controller from \autoref{sec:distr_feedback} are analysed, 
using the standard \eqref{eq:tracking_error} and modified tracking error \eqref{eq:error_dynamic_2_mod}.

\begin{figure}
	\centering
	\begin{subfigure}{\linewidth}
		\centering
		\resizebox {\linewidth} {!} {
			\input{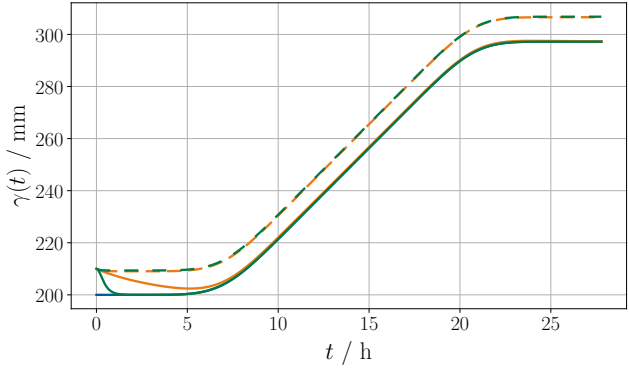}
		}
		\caption{Phase boundary position}
		\label{fig:state_fb_zi}
	\end{subfigure}
	\begin{subfigure}{\linewidth}
		\centering
		\resizebox {\linewidth} {!} {
			\input{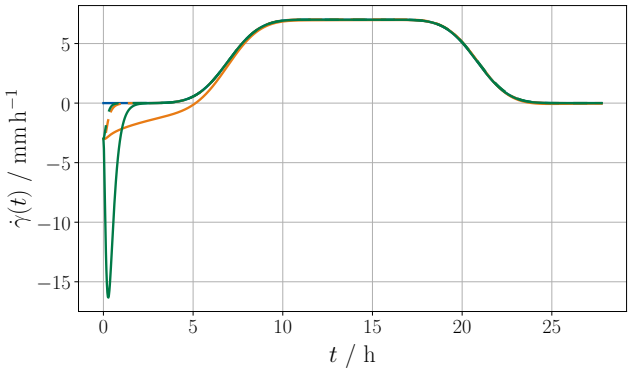}
		}
		\caption{Growth velocity}
		\label{fig:state_fb_vi}
	\end{subfigure}
	\caption{
		Comparison of the tracking behaviour of the reference (blue) between collocated (orange)
		and distributed feedback (green) using the original (solid) or modified error (dashed),
		respectively.
		Note that the orange and green dashed lines are nearly equal.
	}
	\label{fig:state_fb}
\end{figure}

As \autoref{fig:state_fb_zi} shows, the ``fixed'' collocated feedback (orange, solid) successfully
corrects the initial error in the phase boundary position and tracks the reference.
To do this however, a part of the already grown crystal has to be remolten which is to be avoided.
In contrast, the collocated feedback with the shifted error system (orange, dashed) ignores the
error in $\zi$ and makes no attempts on remelting the crystal.
Furthermore, in \autoref{fig:state_fb_vi} it can be seen that this variant corrects the growth rate
error faster than its fixed counterpart. 
However, a drawback that remains for this controller is that due to the simple reference shifting,
a larger crystal is obtained at the end of the process if no further logic is superimposed.

Now to the standard variant of the distributed feedback (green, solid). 
As it can be seen in \autoref{fig:state_fb_zi}, the error in $\zi$ is successfully corrected and
the growth target is reached.
Nevertheless, \autoref{fig:state_fb_vi} shows a severe spike in the growth rate, originating from
the swift correction of $\zi$, thus remelting the crystal.
Opposed to this, the version with the modified tracking error (green, dashed) tolerates the
initial deviation in $\zi$ and continues tracking the trajectory of $\vi$.
However, as for the shifted variant of the collocated feedback, the deviation in $\zi$ still
appears at the end.
The control parameters for these simulations can be found in \ref{app:c}.

\subsection{Observer}
\begin{figure}
	\centering
	\resizebox {\linewidth} {!} {
		\input{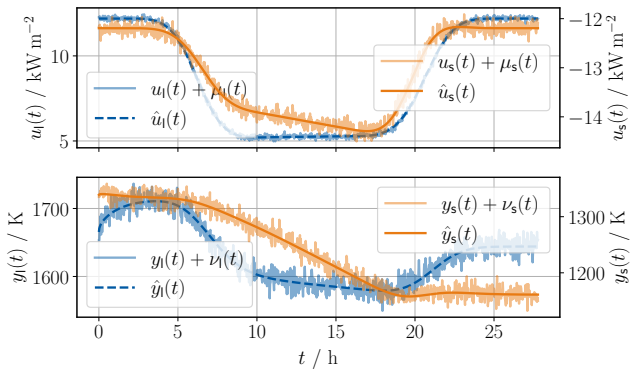}
	}
	\caption{Inputs (top) and outputs (bottom) of the original system (dark) and their disturbed 
		counterparts (light), used by the observer.}
	\label{fig:obs_inputs}
\end{figure}

\begin{figure*}
	\centering
	\begin{subfigure}{.5\linewidth}
		\centering
		\resizebox {\linewidth} {!} {
			\input{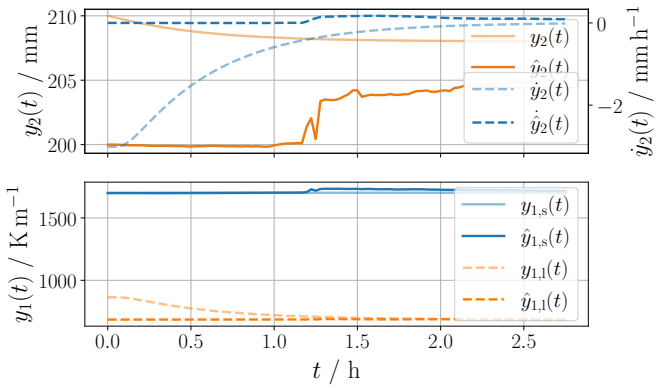}
		}
		\caption{Beginning
		}
		\label{fig:obs_val_red}
	\end{subfigure}%
	\begin{subfigure}{.5\linewidth}
		\centering
		\resizebox {\linewidth} {!} {
			\input{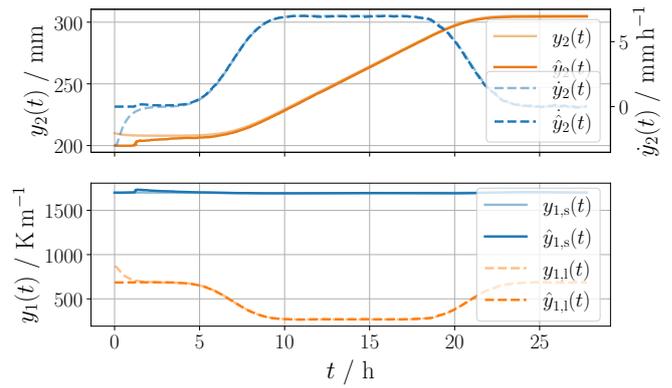}
		}
		\caption{Complete simulation}
		\label{fig:obs_val_full}
	\end{subfigure}%
	\caption{Estimates from the observer (dark), compared to the evolution of the real system's
		variables (light) for the phase boundary and growth rate (top) as well as the gradients
		at the phase boundary for crystal and melt.
	}
	\label{fig:obs_validation}
\end{figure*}

To analyse the observer performance a system under pure feedforward control is considered.
For this case, the initial state estimate $\hat{\bs{\chi}}_N(0)$ of the observer is specified using the
reference trajectory $\bs{y}_{\textup{r}}(0)$. 
To
examine the robustness against disturbances, the real system starts with the initial errors
$\zierro = \SI{100}{\milli\metre}$ and $\vierro = \SI{-3}{\milli\metre\per\hour}$ for both, the crystallisation
interface and the growth rate, respectively.  
Furthermore, the process disturbances $\ud$ and $\yd$ are realised by zero-mean normal distributed noise, distorting the input and output measurements as illustrated in \autoref{fig:obs_inputs}.
As \autoref{fig:obs_val_full} displays, the state estimate quickly converges against the real one
and the system state is properly tracked afterwards.
Due to the different scales of the components in $\statef$, internally a scaled version
$\statefs = \scalemat\statef$ has been used, with $\scalemat$ given in \ref{app:c}.

\subsection{Complete control system}
\begin{figure}
	\centering
	\resizebox {\linewidth} {!} {
		\input{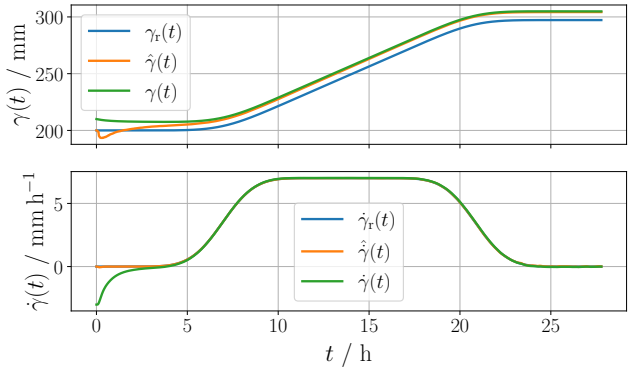}
	}
	\caption{Resulting trajectories of the phase boundary (top) and the growth rate (bottom) for 
	the complete control loop with distributed output feedback (green),
	compared to the reference (blue) and the estimated (orange).
}
	\label{fig:bc_with_ob}
\end{figure}

To examine the performance of the feedback controller when supplied with the state estimates from the 
observer instead of the real system state, a similar setup is used.
Particularly, the estimate is generated by an observer of type \eqref{eq:nonlin_ss_dist} which is
using the model introduced in \autoref{ssec:flat_sys}.
By way of example, the distributed feedback controller with modified tracking error is used to close
the control loop.
As the bottom plot in \autoref{fig:bc_with_ob} shows, the closed loop performs as expected.
%

\section{Conclusion and outlook}

In this contribution, reference tracking control strategies for the \gls{vgf} process,
modelled as a two-phase \gls{sp} have been presented. 
Based on the process demands not to remelt the already solidified domain, two different control 
approaches have been developed.
The performance of components of the control system has been proven by
a simulation study.

A drawback of the proposed solutions is the convergence of the utilised series for smaller transition times.
While this is not a problem for the growth of \gls{gaas}, it may cause complications for the
production of other materials. A direct alternative would be to use so called $(N,\xi)$-approximate $k$-sums
as introduced in \cite{Meurer2004}. 
However, another promising approach is the control via a time-variant backstepping transformation c.f~\cite{Meurer2009} which is currently under investigation and will be covered in a forthcoming publication.


\section*{Acknowledgments}
This work has been funded by the Deutsche Forschungsgemeinschaft (DFG) [project numbers WI 4412/1-1, FR 3671/1-1].

\bibliography{bibliography}

\appendix
\section{Transformations}
%

\subsection{Moving reference system}
\label{app:trafo_mrs}
The coordinate transform
\begin{equation*}
	\ttransz = \torigz \quad \text{with} \quad \zt \coloneqq z - \zi
\end{equation*}
Taking the partial derivative w.r.t.~$z$ gives
\begin{equation}
	\partial_z^2 \torigz = \partial_{\zt}^2 \ttransz
	\label{eq:coord_trans_ddz}
\end{equation}
while according to the chain rule, the time derivative becomes
\begin{align}
	\partial_t \torigz
&= \frac{\mathrm{d}}{\mathrm{d}t} \Big(\ttransz\Big)
	= \partial_{\zt}\ttrans(\zt, t)\partial_t \zt(t) + \partial_t \ttransz \nonumber\\
	&= -\partial_{\zt} \ttransz \vi + \partial_t \ttransz
	\label{eq:coord_trans_dt}
\end{align}
where $\vi$ denotes the growth rate. Inserting \eqref{eq:coord_trans_dt} and
\eqref{eq:coord_trans_ddz} into \eqref{eq:lin_heat_eqs} yields the transformed generic system
\begin{align*}
	\partial_t \ttransz &= \hd\partial_{\zt}^2\ttransz + \vi \partial_{\zt} \ttransz \\
	\hc\partial_{\zt} \ttrans(\zb - \zt, t) &= \delta\inp \\
	\ttrans(0, t) &= \tm \:.
\end{align*}%

\subsection{Shifted error system}
\label{app:trafo_err}
The evolution of the ``shifted'' error can be described by taking the time derivative of 
\eqref{eq:shifted_error} for the real and the planned profile
\begin{align*}
	\partial_t \terrz 
	&= \partial_t \torigz - \frac{\mathrm{d}}{\mathrm{d}t}\big(\treft\big) \\
	&= \partial_t \torigz 
		- \partial_t \treft\\
	&\quad + \partial_z\treft \dvi
	\intertext{and then substituting the \glspl{pde} according to \eqref{eq:lin_heat_eqs},
		which hold for the real as well as for the reference system:
	}
	\partial_t \terrz 
	&= \hd\partial_z^2 \torigz 
		- \hd\partial_z^2 \treft\\
	&\quad + \dvi \partial_z \treft \:.
	\intertext{Making use of definition \eqref{eq:shifted_error}, one arrives at}
	\partial_t \terrz &= \hd\partial_z^2 \terrz 
		+ \dvi \partial_z \treft \:,
	\intertext{while the related boundary conditions are given by}
	\partial_z \terr(\zb, t) 
		&= \partial_z\torig(\zb, t) -  \partial_z\tref(\zb-\dzi, t)\\
		&= \frac{\delta}{\hc} u(t) - \partial_z \tref(\zb - \dzi) \\
	\terr(\zi, t) &= \torig(\zi, t) -  \tref(\zi-\dzi, t)\\
		&= T_m - T_m = 0 \:.
\end{align*}
Hence, the resulting error system is governed by
\begin{align*}
	\partial_t \terrz &= \hd\partial_z^2 \terrz +  \dvi \partial_z\treft\\
	\partial_z \terr(\zb, t) &= \frac{\delta}{\hc} u(t)
	- \partial_z \tref(\zb - \dzi, t) \\
	\terr(\zi, t) &= 0 \:.
\end{align*}%

\section{Stability analysis}
\label{app:a}

Firstly, $\Vo$ is decomposed into
\begin{equation}
	\Vo =\Vs + \Vl
	\label{eq:comp_lyap_func}
\end{equation}
with $\Vs = \frac{1}{2}\int\limits_{\zbs}^{\zi} \terrszsq \, \textup{d}z$ and 
$\Vl = \frac{1}{2}\int\limits_{\zi}^{\zbl} \terrlzsq \, \textup{d}z$.
For the sake of brevity, the next steps will focus on the generic function $\Vg$ since they are
similar for $\Vs$ and $\Vl$. Differentiation of 
\begin{equation}
	\Vg = \frac{\bfd}{2}\int\limits_{\zi}^{\zb} \terrgzsq \, \textup{d}z
	\label{eq:gen_ljap_func}
\end{equation}
leads to
\begin{equation*}
	\Vgt 
	= -\frac{\bfd}{2}\vi\terrgsq(\zi, t)
	+ \bfd\int\limits_{\zi}^{\zb} \terrgz \partial_t \terrgz \, \textup{d}z \:.
\end{equation*}
Using \eqref{eq:error_system_bc2} and substituting the system dynamics
\eqref{eq:error_system_pde} one obtains
\begin{align*}
	\Vgt = \bfd\int\limits_{\zi}^{\zb} \terrgz &\Big(\hd\partial_z^2 \terrz  \\
	&+ \dvi\partial_z\trefgt \Big) \, \textup{d}z \:.
\end{align*}
Integration by parts of the first summand yields
\begin{align}
	\begin{split}
		\Vgt &= 
		\hd\bfd\big[\terrgz\partial_z\terrgz\big]_{\zi}^{\zbs} \\
		&\hphantom{=} -\hd\bfd\int\limits_{\zi}^{\zbs} \left(\partial_z \terrgz\right)^2 \, \textup{d}z\\
		&\hphantom{=}  +\bfd\dvi\int\limits_{\zi}^{\zbs} 
			\terrgz \partial_z\trefgt \,\textup{d}z \:.\nonumber
	\end{split}
	\intertext{Using the boundary conditions \eqref{eq:error_system_bc1}, 
		\eqref{eq:error_system_bc2} as well as the feedback law \eqref{eq:bc_feedback} gives}
	\begin{split}
		&= -\hd\cgaing\terrgzbsq
		- \hd\bfd\int\limits_{\zi}^{\zb} \left(\partial_z \terrg\right)^2 \, \textup{d}z \\
		&\quad+\bfd\dvi\int\limits_{\zi}^{\zb} \terrg \partial_z\trefgt \, \textup{d}z 
	\end{split}
	\label{eq:vs_dt}
\end{align}
since $\bfd^2 = 1$.
To reassemble $\Vgt$ in the expression, the first term has to be rearranged. 
Therefore, a slightly modified version of the \gls{pi} is introduced, based on the one
given in \cite{book:KrsticBackstepping}.

\paragraph{Poincaré inequality}
Consider the partially integrated term
\begin{align*}
    \int\limits_{\zi}^{\zb}\terrgzsq \,\textup{d}z 
    &= \big[z \terrgzsq \big]_{\zi}^{\zb}
       - 2\int\limits_{\zi}^{\zb}  z \terrgz \partial_z\terrgz \, \textup{d}z \\
    &= \zb\terrgzbsq 
	- 2\int\limits_{\zi}^{\zb}  z \terrgz \partial_z\terrgz \, \textup{d}z 
	\intertext{which, by making use of the \gls{csi}, can be estimated as}
    \int\limits_{\zi}^{\zb}\terrgzsq \,\textup{d}z 
    &\le 
	\overbrace{
			\sqrt{\int\limits_{\zi}^{\zb} \terrgzsq \, \textup{d}z }
		}^{\eqqcolon a}
        \overbrace{
			2\sqrt{\int\limits_{\zi}^{\zb} \left(z \partial_z\terrgz\right)^2 \, \textup{d}z}
		}^{\eqqcolon b} \\
	&\quad +\zb\terrgzbsq\:.
		\intertext{Using \gls{yi} $ab \le \frac{1}{2\sigma}a^2 + \frac{\sigma}{2}b^2$
	with $\sigma = 1$ one arrives at}
    \int\limits_{\zi}^{\zb}\terrgzsq \,\textup{d}z 
    &\le \frac{1}{2}\int\limits_{\zi}^{\zb} \terrgzsq \, \textup{d}z
        + 2\int\limits_{\zi}^{\zb} \left(z \partial_z\terrgz\right)^2 \, \textup{d}z \\
		&\quad +\zb\terrgzbsq \:,
	\intertext{which after rearranging can be further estimated by}
	\frac{1}{2}\int\limits_{\zi}^{\zb}\terrgz \,\textup{d}z 
    &\le 2\int\limits_{\zi}^{\zb} \left(z \partial_z\terrgz\right)^2 \, \textup{d}z
		+\zb\terrgzbsq \\
	&\le 2\zbl^2 \int\limits_{\zi}^{\zb} \left(\partial_z\terrgz\right)^2 \, \textup{d}z
		+\zb\terrgzbsq \:,
\end{align*}
since $\zbs \le \zi \le \zbl \forall t$.
Further rearrangement finally provides the required inequality:
\begin{equation}
    -\int\limits_{\zi}^{\zb} \left(\partial_z\terrgz\right)^2 \, \textup{d}z
	\le -\frac{1}{4\zbl^2}\int\limits_{\zi}^{\zb}\terrgzsq \,\textup{d}z
	-\frac{\zb}{2\zbl^2}\terrgzbsq \:.
    \label{eq:pc_ineq_s}
\end{equation}

\paragraph{Reassembly}
Substitution of \eqref{eq:pc_ineq_s} in \eqref{eq:vs_dt} yields
\begin{align*}
	\Vgt 
	&\le 
	-\hd\terrgzbsq\left(\cgaing + \frac{\bfd\zb}{2\zbl^2}\right)
	-\frac{\hd\bfd}{4\zbl^2}\int\limits_{\zi}^{\zb}\terrgzsq \,\textup{d}z  \\
	&\hphantom{\le} 
	+ \bfd\dvi\int\limits_{\zi}^{\zb} \terrgz \partial_z\trefgt \, \textup{d}z \:.
	\intertext{Thus, comparison with \eqref{eq:gen_ljap_func} leads to}
	&\le 
	-b_{\circ}\terrgzbsq
		-\frac{\hd}{2\zbl^2} \Vg   \\
	&\hphantom{\le} 
	+	\dvi\int\limits_{\zi}^{\zb} \terrgz \partial_z\trefgt \, \textup{d}z
\end{align*}
with the positive constant $b_{\circ} = \hd\left(\cgaing + \frac{\bfd\zb}{2\zbl^2}\right)$
by an appropriate choice of $\cgaing$.
Hence, by substituting the results for $\Vs$ and $\Vl$ in \eqref{eq:comp_lyap_func}, $\Vot$
can be expressed as 
\begin{align}
	\begin{split}
		\Vot
		&\le 
		-b_{\textup{s}}\terrszbsq
		-b_{\textup{l}}\terrlzbsq
		-C \Vo \\
		&\hphantom{\le}+	\dvi\int\limits_{\zbs}^{\zbl} \terrz \partial_z\treft \, \textup{d}z%
	\end{split}%
	\label{eq:ljap_cand_dt_5}%
\end{align}%
with $C = \frac{1}{2\zbl^2}\max\left\{\hds, \hdl\right\} \ge 0$.

\paragraph{Simplified variant}
Choosing the reference profile to be constant w.r.t.~the spatial dimension $z$ makes the integral
term in \eqref{eq:ljap_cand_dt_5} vanish. This can be achieved by using the reference 
$\tref^0 \equiv \tm$, which may be the case if the process should be brought to halt. 
The resulting derivative
\begin{equation}
	\Vot
	\le 
	-b_{\textup{s}}\terrszbsq
	-b_{\textup{l}}\terrlzbsq
	-C \Vo
	\label{eq:ljap_cand_dt_5_simp}
\end{equation}
is obviously negative definite, yielding uniformly, globally, asymptotic stability 
of the system \eqref{eq:error_system_pde} with respect to $\targetset$ for this case.

\paragraph{General variant}
Ignoring the boundary terms in \eqref{eq:ljap_cand_dt_5} and focussing on the last term,
by using \gls{csi}, the estimation
\begin{equation*}
	\dvi\int\limits_{\zbs}^{\zbl} \terrz \partial_z\treft \, \textup{d}z 
	\le \left|\dvi\right|
		K
		\sqrt{\Vo}
\end{equation*}
is obtained, where 
$K = \sqrt{2} \left(\zbl-\zbs\right) \max_{z, t} \left|\partial_z\trefz \right|$.
Summarizing, $\Vot)$ is bounded by
\begin{equation}
	\Vot \le -C \Vo + \left| \dvi \right| K \sqrt{\Vo} \:.
	\label{eq:ljap_cand_dt_6}
\end{equation}
Furthermore, for every $\Vo$ a scaling $\nu > 0$ can be found such that
\begin{equation*}
	\Vo \ge \nu\sqrt{\Vo}
\end{equation*}
holds. 
As a consequence, for $\Vo \ge \nu^2$ the estimate
\begin{equation}
	\Vot \le \left(\left| \dvi \right| \bar K  -C\right)\Vo
	\label{eq:ljap_cand_dt_7}
\end{equation}
with $\bar K = \frac{K}{\nu}$ can be used.
Finally by using Gronwalls lemma, a growth bound for $\Vo$ is given by
\begin{equation}
	\Vo \le V(0) \exp(\bar K \Vdzi -Ct) \:.
	\label{eq:ljap_cand_sol}
\end{equation}
Thus, \eqref{eq:ljap_cand} is decreasing for $\Vo \ge \nu^2$ if 
\begin{equation}
	\Vdzi \le \frac{C}{\bar K}t
	\label{eq:variation_bound}
\end{equation}
holds, hence, if there exists an upper bound for the total variation $\Vdzi$ that grows linear in time.

\section{Matrices and Vectors}
\label{app:b}

\begin{align*}
	\bs{P}_{0,0} &= \begin{psmallmatrix}
		\inner{\bar\varphi_0(\bar z) }{\bar \varphi_{0}(\bar z)} 
		&\dotsc 
		&\inner{\bar\varphi_{N-1}(\bar z) }{\bar \varphi_0(\bar z)}\\
		\vdots &\ddots &\vdots \\
		\inner{\bar\varphi_0(\bar z) }{\bar \varphi_{N-1}(\bar z)} 
		&\dotsc 
		&\inner{\bar\varphi_{N-1}(\bar z) }{\bar \varphi_{N-1}(\bar z)}
	\end{psmallmatrix} \\
	\bs{P}_{1,0} &= \begin{psmallmatrix}
		\inner{\bar z\partial_{\bar z}\bar\varphi_0(\bar z) }{\bar \varphi_{0}(\bar z)} 
		&\dotsc 
		&\inner{\bar z\partial_{\bar z}\bar\varphi_{N-1}(\bar z) }{\bar \varphi_0(\bar z)}\\
		\vdots &\ddots &\vdots \\
		\inner{\bar z\partial_{\bar z}\bar\varphi_0(\bar z) }{\bar \varphi_{N-1}(\bar z)} 
		&\dotsc 
		&\inner{\bar z\partial_{\bar z}\bar\varphi_{N-1}(\bar z) }{\bar \varphi_{N-1}(\bar z)}
	\end{psmallmatrix} \\
	\bs{P}_{1,1} &= \begin{psmallmatrix}
		\inner{\partial_{\bar z}\bar\varphi_0(\bar z) }{\bar \varphi_{0}(\bar z)} 
		&\dotsc 
		&\inner{\partial_{\bar z}\bar\varphi_{N-1}(\bar z) }{\bar \varphi_0(\bar z)}\\
		\vdots &\ddots &\vdots \\
		\inner{\partial_{\bar z}\bar\varphi_0(\bar z) }{\bar \varphi_{N-1}(\bar z)} 
		&\dotsc 
		&\inner{\partial_{\bar z}\bar\varphi_{N-1}(\bar z) }{\bar \varphi_{N-1}(\bar z)}
	\end{psmallmatrix} \\
	\bs{q}_{1,0} &= \begin{psmallmatrix}
		\inner{\bar z\partial_{\bar z}\bar\varphi_N(\bar z) }{\bar \varphi_{0}(\bar z)} \\
		\vdots \\
		\inner{\bar z\partial_{\bar z}\bar\varphi_{N}(\bar z) }{\bar \varphi_{N-1}(\bar z)}
	\end{psmallmatrix},
	\bs{q}_{1,1} = \begin{psmallmatrix}
		\inner{\partial_{\bar z}\bar\varphi_N(\bar z) }{\bar \varphi_{0}(\bar z)} \\
		\vdots \\
		\inner{\partial_{\bar z}\bar\varphi_{N}(\bar z) }{\bar \varphi_{N-1}(\bar z)}
	\end{psmallmatrix} \\
\end{align*}

\section{Control parameters}
\label{app:c}
The parameters of the collocated feedback have been chosen as $\cgains = \cgainl = \SI{20}{\per\metre}$,
while the distributed feedback was parametrised with 
$\kappa_{1,0}=\SI{2e-6}{\per\square\metre}$,
$\kappa_{1,1}=\SI{3e-3}{\per\metre}$ as well as 
$\kappa_{2,0}=\SI{6e-9}{\per\cubic\second}$,
$\kappa_{2,1}=\SI{1.1e-5}{\per\square\second}$,
$\kappa_{2,2}=\SI{6e-3}{\per\second}$ for the original and
$\tilde\kappa_{2,0}=\SI{6e-6}{\per\square\second}$,
$\tilde\kappa_{2,1}=\SI{5e-3}{\per\second}$
for the modified error.

The scaling matrix for $\statefs[5]$ was chosen as 
\begin{equation*}
	\scalemat[5] = \diag{\SI{1e-3}{}, \SI{1e-3}{\metre}, \SI{1e7}{}, \SI{1e10}{\metre}, \SI{1e13}{\square\metre}}.
\end{equation*}

\end{document}